\documentclass[11pt]{article}
\usepackage{amssymb}
\usepackage{amsmath}
\usepackage{graphicx}
\usepackage{lipsum}

\title{ }
\author{}
\date{}

\newtheorem{theorem}{Theorem}[section]

\newtheorem{lemma}[theorem]{Lemma}

\newtheorem{coro}[theorem]{Corollary}

\newcommand{\qed}{\hspace*{\fill} \rule{7pt}{7pt}}

\newcommand{\E}{{\mathbb{E}}}
\newcommand{\F}{{\mathbb{F}}}

\newcommand{\X}{{\mathbf{X}}}
\newcommand{\pr}{\mathbf{Pr}}

\newcommand{\sP}{\mathcal{P}}

\newcommand\blfootnote[1]{%
  \begingroup
  \renewcommand\thefootnote{}\footnote{#1}%
  \addtocounter{footnote}{-1}%
  \endgroup
}

\begin{document}

\title{Some extremal results on complete degenerate hypergraphs}
\author{
Jie Ma\footnote{E-mail: jiema@ustc.edu.cn. Research partially supported by NSFC projects 11501539 and 11622110.}
~~~~~
Xiaofan Yuan\footnote{E-mail: yxf518@mail.ustc.edu.cn.}
~~~~~
Mingwei Zhang\footnote{E-mail: zmw309@mail.ustc.edu.cn.}
\\
\footnotesize{School of Mathematical Sciences, University of Science and Technology of China}\\
\footnotesize{Hefei, Anhui 230026, P.R. China.}
}
\date{}

\maketitle
\begin{abstract}
Let $K^{(r)}_{s_1,s_2,\cdots,s_r}$ be the complete $r$-partite $r$-uniform hypergraph and $ex(n, K^{(r)}_{s_1,s_2,\cdots,s_r})$ be
the maximum number of edges in any $n$-vertex $K^{(r)}_{s_1,s_2,\cdots,s_r}$-free $r$-uniform hypergraph.
It is well-known in the graph case \cite{KST,KRS} that $ex(n,K_{s,t})=\Theta(n^{2-1/s})$ when $t$ is sufficiently larger than $s$.
In this note, we generalize the above to hypergraphs by showing that if $s_r$ is sufficiently larger than $s_1,s_2,\cdots,s_{r-1}$
then $$ex(n, K^{(r)}_{s_1,s_2,\cdots,s_r})=\Theta\left(n^{r-\frac{1}{s_1s_2\cdots s_{r-1}}}\right).$$
This follows from a more general Tur\'an type result we establish in hypergraphs,
which also improves and generalizes some recent results of Alon and Shikhelman \cite{AS}.
The lower bounds of our results are obtained by the powerful random algebraic method of Bukh \cite{B}.
Another new, perhaps unsurprising insight which we provide here is that one can also use the random algebraic method to construct
non-degenerate (hyper-)graphs for various Tur\'an type problems.

The asymptotics for $ex(n, K^{(r)}_{s_1,s_2,\cdots,s_r})$ is also proved by Verstra\"ete \cite{Ver} independently with a different approach.
\blfootnote{2010 Mathematics Subject Classification: 05C35, 05C65, 05D99}
\blfootnote{Keywords: Tur\'an number, degenerate hypergraph, random algebraic method, Alon-Shikhelman}
\end{abstract}

\section{Introduction}
Let $r\ge 2$ be an integer and $T,H$ be two $r$-uniform hypergraphs.
An $r$-uniform hypergraph is called {\it $H$-free} if it contains no copy of $H$ as its subhypergraph.
For any integer $n$, let $ex(n,T,H)$ be the maximum number of copies of $T$ in an $n$-vertex $H$-free $r$-uniform hypergraph. In case $T$ is a single edge, this function then converts to the {\it Tur\'an number} $ex(n,H)$ of the hypergraph $H$.

The study of Tur\'an numbers $ex(n,H)$ is the main focus of the extremal graph theory.
This was initiated by the celebrated theorems of Mantel \cite{Mantel} and Tur\'an \cite{Tur},
which determine the precise value of $ex(n,H)$ when $H$ is a complete graph.
For general graphs $H$, Erd\H{o}s-Stone and Simonovits \cite{ESt,ESi,Sim} resolved Tur\'an numbers $ex(n,H)$ asymptotically,
except for bipartite graphs $H$ (which is often called the {\it degenerate graphs}).
Even to date there are still few degenerate graphs the asymptotics of whose Tur\'an numbers are known.
One of such examples is the complete bipartite graph $K_{s,t}$.
A well-known theorem of K\H{o}v\'ari, S\'os and Tur\'an \cite{KST} shows that $ex(n,K_{s,t})=O(n^{2-1/s})$ for any integers $t\ge s$.
For $s=2, 3$, matched lower bounds were found in \cite{ERS,Bro} respectively;
for other values of $s$, this bound was known to be tight when $t$ is sufficiently larger than $s$,
which was first proved by Koll\'ar, R\'onyai and Szab\'o \cite{KRS}
and then slightly improved to $t>(s-1)!$ by Alon, R\'onyai and Szab\'o \cite{ARS}.
Recently, Blagojevi\'c, Bukh and Karasev \cite{BBK} and Bukh \cite{B} used the random algebraic method
to give different constructions which yield the same lower bound $ex(n,K_{s,t})=\Omega(n^{2-1/s})$ as in \cite{KRS,ARS},
provided that $t$ is sufficiently large. For more extremal results on generate graphs, we refer to the survey \cite{FS}.

For the function $ex(n,T,H)$ where $T, H$ are graphs and $T$ is not an edge,
there are only sporadic results such as \cite{Erd,Boll,GL} in the literature,
until recently Alon and Shikhelman \cite{AS} systematically investigate this general function
and obtain a number of results on complete graphs, complete bipartite graphs and trees.
Among other results, they \cite{AS} proved that for $s\ge 2m-2$ and $t\ge (s-1)!+1$,
\begin{align}\label{equ:AS-Km}
ex(n,K_m,K_{s,t})=\Theta\left(n^{m-\frac{m(m-1)}{2s}}\right)
\end{align}
and for $(a-1)!+1\le b<(s+1)/2$ and $t\ge s$,
\begin{align}\label{equ:AS-Kab}
ex(n,K_{a,b},K_{s,t})=\Theta\left(n^{a+b-ab/s}\right).
\end{align}

By contrast with the graph case, very little was known for the hypergraph Tur\'an problems.
(For instance, it is not known the value of $ex(n,H)$ when $H$ is a complete
$r$-uniform hypergraph on $t$ vertices for {\it any} pair of integers $t>r>2$, even asymptotically.)
An $r$-uniform hypergraph $H$ is called {\it degenerate} if it is {\it $r$-partite},
i.e., its vertices can be partitioned into $r$ parts such that each edge contains exactly one vertex from each part.
Let $K^{(r)}_{s_1,s_2,...,s_r}$ be the complete $r$-partite $r$-uniform hypergraph with parts of sizes $s_1,...,s_r$ respectively.
A classical result of Erd\H{o}s \cite{Er64} indicates that $ex(n, K^{(r)}_{s_1,s_2,...,s_r})\le O\left(n^{r-1/(s_1\cdots s_{r-1})}\right)$ for all $s_1\le s_2\le \cdots \le s_r$.
It was probably widely believed but not stated formally until in \cite{Mu}, where Mubayi conjectured that for all $s_1\le s_2\le \cdots \le s_r$,
$$ex(n, K^{(r)}_{s_1,s_2,...,s_r})=\Theta\left(n^{r-1/(s_1\cdots s_{r-1})}\right);$$
some sharp constructions were also given in \cite{Mu} in case $s_1=\cdots=s_{r-2}=1$.
The special case $ex(n, K^{(r)}_{2,2,...,2})$ intrigues many researchers and is often referred as {\it the box problem}.
In \cite{GRS} Gunderson, R\"odl and Sidorenko provided a construction for $ex(n, K^{(r)}_{2,2,...,2})$ which beats the canonical lower bound from the probabilistic deletion method.
As for $r=3$, the current best bounds are $\Omega(n^{8/3})\le ex(n, K^{(3)}_{2,2,2})\le O(n^{11/4})$, where the lower bound was obtained by Katz, Krop and Maggioni in \cite{KKM}.

Our first result generalizes \eqref{equ:AS-Kab} to hypergraphs.
\begin{theorem} \label{thm:Ksssp}
  For any integers $a_1,a_2,...,a_r, s_1, s_2,..., s_{r-1}$ satisfying that $a_1<s_1\le s_2\le \cdots\le s_{r-1}$
  and $a_i\le s_{i-1}$ for each $2\le i\le r$, there exists a constant $c>0$ such that the following holds. For any integer $p\ge c$, we have
  $$ex(n,K^{(r)}_{a_1,a_2,...,a_r},K^{(r)}_{s_1,s_2,...,s_{r-1},p})=\Theta \left(n^{(a_1+a_2+\cdots +a_r)-\frac{a_1 a_2 \cdots a_r}{s_1s_2\cdots s_{r-1}}}\right).$$
\end{theorem}

\noindent
From this, we can promptly obtain the following asymptotic bound for Tur\'an numbers of complete degenerate hypergraphs\footnote{See the remark in the end of Section 3.}, which partially confirms the above conjecture of Mubayi.
\begin{coro}\label{cor:Ksssp}
  For any positive integers $s_1, s_2,..., s_{r-1}$, there exists a constant $c>0$ such that for any integer $p\ge c$,
  $$ex(n,K^{(r)}_{s_1,s_2,...,s_{r-1},p})=\Theta \left(n^{r-\frac{1}{s_1s_2\cdots s_{r-1}}}\right).$$
\end{coro}

\noindent This result is also obtained by Verstra\"ete \cite{Ver} independently.

Also when restricting to $r=2$, Theorem \ref{thm:Ksssp} improves \eqref{equ:AS-Kab}
by weakening the relations between $a,b$ and $s$, at the cost of requiring $t$ to be even larger.

\begin{coro}\label{cor:Kab-Kst}
For any integers $a<s$ and $b\le s$, there exists a constant $f(a,b,s)>0$ such that for any integer $t\ge f(a,b,s)$,
we have $ex(n,K_{a,b},K_{s,t})=\Theta\left(n^{a+b-ab/s}\right).$
\end{coro}

Our construction of the lower bound in Theorem \ref{thm:Ksssp} is from the random algebraic method,
which was initialed in \cite{BBK,B} and developed in \cite{C,BC} quite recently.
Similarly as in \cite{BBK,B,C,BC}, it also suffices for us to construct a degenerate hypergraph,
for the purpose of counting complete degenerate hypergraphs.
However, if one wants to count the number of copies for some non-degenerate hypergraph $H$,
then it will be necessary to consider non-degenerate examples.
Using a variant of the random algebraic method, we show the following result for general $H$.

\begin{theorem}\label{thm:lower}
For any integers $s_1,s_2,...,s_{r-1}$ and any $r$-uniform hypergraph $H$,
there exists a constant $c>0$ such that for any integer $p\ge c$,
$$ex(n,H,K^{(r)}_{s_1,s_2,...,s_{r-1},p})=\Omega \left(n^{v-\frac{e}{s_1s_2\cdots s_{r-1}}}\right),$$
where $v=|V(H)|$ and $e=|E(H)|$.
\end{theorem}

\noindent In the graph case, if we choose $H$ to be a complete graph $K_m$, then together with an upper bound from \cite{AS} that $ex(n,K_m,K_{s,t})=O\left(n^{m-\frac{m(m-1)}{2s}}\right)$ for $t\ge s\ge m-1$,
we can deduce the following analog of \eqref{equ:AS-Km}, by relaxing the condition to $s\ge m-1$.

\begin{coro}\label{cor:Km-Kst}
For any integers $s, m$ with $s+1\ge m\ge 2$, there exists a constant $f(s,m)>0$ such that for any integer
$t\ge f(s,m)$, we have
$ex(n,K_m,K_{s,t})=\Theta\left(n^{m-\frac{m(m-1)}{2s}}\right).$
\end{coro}

The rest of this paper is organized as follows. In Section 2,
we show the lower bounds of our results by establishing Theorem \ref{thm:lower}.
In Section 3, we prove the upper bound of Theorem \ref{thm:Ksssp} for degenerate hypergraphs.
In Section 4, we conclude by some remarks and a related open problem.
For a set $V$ and an integer $r$, we write $[r]:=\{1,2,\cdots, r\}$ and $\binom{V}{r}$
as the family of all subsets of size $r$ in $V$.

\section{Lower bounds}
In this section, we prove Theorem \ref{thm:lower} by providing a random algebraic construction.
The main part of the proof will follow the line of \cite{B,C,BC} closely.
For given positive integers $s_1, s_2, \cdots, s_{r-1}$ and an $r$-uninform hypergraph $H$, throughout this section let $$v=|V(H)|,~~ e=|E(H)|, ~~b=\prod_{i=1}^{r-1}s_i ~~ \text{ and } ~~s=b\left(\sum_{i=1}^{r-1}s_i-1\right)+e+1.$$
Let $q$ be a sufficiently large prime power (compared to the above parameters) and let $\F_q$ be the finite field of order $q$.

Let $\X^i=(X^i_1,X^i_2,\cdots,X^i_b)\in \F_q^b$ for each $i\in [r]$.
Consider polynomials $f\in \F_q[\X^1,\X^2,...,\X^r]$ with $rb$ variables over $\F_q$.
We say such a polynomial $f$ has degree at most $d$ in $\X^i$, if each of its monomials with respect to $\X^i$, say $(X^i_1)^{\alpha_1}(X^i_2)^{\alpha_2}\cdots (X^i_b)^{\alpha_b}$, satisfies $\sum_{j=1}^b \alpha_j\le d$.
And a polynomial $f$ is called {\it symmetric}, if exchanging $\X^i$ with $\X^j$ for every $1\le i<j\le r$ will not affect the value of $f$.
It will be convenient to view the domain of symmetric polynomials as the family $\binom{\F_q^b}{r}$.
Given such a symmetric polynomial $f$, we then can define an $r$-uniform hypergraph $G_f=(V,E)$ as following:
the vertex set $V$ is a copy of $\F_q^b$, and every $\{u^1,\cdots,u^r\}\in \binom{V}{r}$ is an edge of $G_f$ if and only if $f(u^1,\cdots,u^r)=0$.

Let $\sP\subset \F_q[\X^1,\X^2,...,\X^r]$ be the set of all symmetric polynomials of degree at most $bs$ in $\X^i$ for every $1\le i\le r$.
Choose a polynomial $f$ from $\sP$ uniformly at random and let $G=G_f$ be the associated $r$-uniform hypergraph.
Our goal is to show that averagely this $G$ contains many copies of $H$ but very few copies of $K^{(r)}_{s_1,s_2,...,s_{r-1},p}$, assuming $p$ is sufficiently large;
then one can use the deletion method to obtain a subhypergraph of $G$ which is $K^{(r)}_{s_1,s_2,...,s_{r-1},p}$-free and yet has expected number of copies of $H$.

We shall first present some lemmas which are needed later in the proof.
A set of points in $\F_q^b$ is called {\it simple} if the first coordinates of all the points are distinct.
The analogous proof of the coming lemma can be found in \cite{B}.

\begin{lemma}\label{lem:simple}
Let $V\subset\F_q^b$ be a set of points. Suppose $\binom{|V|}{2} < q$. Then there exists an invertible linear transformation $T:\F_q^b\to \F_q^b$ such that $TV\subset\F_q^b$ is simple.
\end{lemma}

\begin{pf}
Let $u=(u_1,\ldots,u_b)$ and $v=(v_1,\ldots,v_b)\in\F_q^b$ be two points. We write $u\cdot v$ for $\sum u_i v_i\in\F_q,$
and for each $u\in \F_q^b$, define the linear function $L_u: \F_q^b\to \F_q$ by assigning $L_u(x)=u\cdot x$ for every $x\in \F_q^b$.
It suffices to find a point $u\in\F_q^b$ such that $L_u(v-v')\neq 0$ for every $\{v,v'\}\in \binom{V}{2}$.
That is, $u\notin\bigcup\ker L_{v-v'}$ over all pairs $\{v,v'\}\in \binom{V}{2}$.
Since the kernel of each linear function $L_{v-v'}$ contains $q^{b-1}$ points, the existence of the point $u$ then follows that
$|\bigcup\ker L_{v-v'}|\leq \sum|\ker L_{v-v'}|=\binom{|V|}{2} q^{b-1}<q^b.$\qed
\end{pf}

\medskip

The following lemma will be crucial, which says that for any set $U\subseteq \binom{\F_q^b}{r}$ of bounded size,
the probability that $U$ appears in the kernel of a randomly chosen polynomial $f\in \sP$ is precisely $1/q^{|U|}$.
In this sense, we see that $f$ indeed behaves randomly.

\begin{lemma}\label{lem:prob}
Given a set $U\subseteq \binom{\F_q^b}{r}$, let $V\subset\F_q^b$ be the set consisting of all points appeared as an element of an $r$-tuple in $U$.
Suppose that $\binom{|U|}{2} < q, ~\binom{|V|}{2} < q$ and $|U|\le bs$. If $f$ is a random polynomial chosen from $\sP$, then
$$\pr[f(u^1,u^2,...,u^r)=0,\ \forall\{u^1,u^2,...,u^r\}\in U]=q^{-|U|}.$$
\end{lemma}

\begin{pf}
By Lemma \ref{lem:simple}, there exists an invertible linear transformation $T:\F_q^b\to \F_q^b$ such that $TV\subset\F_q^b$ is simple.
Then $T$ induces an invertible linear transformation $T^\ast:\sP\to\sP$ by letting $$f(\X^1,\X^2,...,\X^r)\mapsto f(T\X^1,T\X^2,...,T\X^r), \text{ for each } f\in \sP.$$
Therefore, it will suffice for us to consider that $V$ is a simple set.

Now observe that any $\{u^1,u^2,...,u^r\}\in U$ is uniquely determined by $\{u^1_1,u^2_1,...,u^r_1\}$, which is uniquely determined by
$$[u^1,u^2,...,u^r]:=\left(\sum_{i=1}^r u^i_1,\ \sum_{1\le i< j\le r}u^i_1u^j_1,\ \ldots\ ,\ \sum_{1\le i_1<i_2<\cdots<i_{r-1}\le r}u^{i_1}_1u^{i_2}_1\cdots u^{i_{r-1}}_1,\ \prod_{i=1}^ru^i_1\right)\in\F_q^r.$$
Since $\binom{|U|}{2}< p$, applying Lemma \ref{lem:simple} (and its proof) to the set consisting of all $[u^1,u^2,...,u^r]$,
we can find an injective linear function $\phi: U\to \F_q$ such that
$$\phi(\{u^1,u^2,...,u^r\})=a_1\cdot\left(\sum_{i=1}^r u^i_1\right)+a_2\cdot\left(\sum_{1\le i< j\le r}u^i_1u^j_1\right)+\cdots +a_r \cdot\prod_{i=1}^ru^i_1$$
for some $a_1, a_2,\cdots, a_r\in \F_q$. Let
$$t=a_1\cdot\left(\sum_{i=1}^r X^i_1\right)+a_2\cdot\left(\sum_{1\le i< j\le r}X^i_1X^j_1\right)+\cdots +a_r \cdot\prod_{i=1}^rX^i_1$$
and $T$ be the linear subspace generated by $1,t,t^2,\ldots t^{|U|-1}$ over $\F_q$.
So $t$ is a symmetric polynomial in $\sP$ of degree at most one in each $\X^i$.
Also as $|U|\le bs$, any polynomial in $T$ has degree at most $bs$ in each $\X^i$.
This shows that $T\subset\sP$, where $|T|=q^{|U|}$. We fix a supplement subspace $W$ of $T$ in $\sP$.

We then decompose $f$ as $f=g+w$, where $g\in T$ and $w\in W$.
Clearly, one can sample $f\in \sP$ uniformly by first sampling $w\in W$ and then sampling $g\in T$.
Hence, to study the system $f(u^1,u^2,...,u^r)=0$ for all $\{u^1,u^2,...,u^r\}\in U$,
it suffices to consider the system of linear equations
$$g(u^1,u^2,...,u^r)=-w(u^1,u^2,...,u^r) \text{ for all } \{u^1,u^2,...,u^r\}\in U,$$
where we view $w$ as a given function and view $g$ as unknown.
In fact, $g$ can be viewed as a polynomial in $\F_q[t]$ of degree at most $|U|-1$ satisfying
$$g(t)=-w(u^1,u^2,...,u^r), \text{ where } t=\phi(u^1,u^2,...,u^r)  \text{ for all } \{u^1,u^2,...,u^r\}\in U.$$
Recall that $\phi$ is injective. Therefore, by the Lagrange Interpolation Theorem,
for every given $w\in W$, there exists a unique solution $g\in T$ to the above system, while there are exactly $q^{|U|}$ polynomials in $T$.
This shows that the probability that $f(u^1,u^2,...,u^r)=0$ for all $\{u^1,u^2,...,u^r\}\in U$ is $\frac{1}{|T|}=q^{-|U|},$ finishing the proof.
\qed
\end{pf}

\medskip

The following lemma of \cite{B} is also important. It indicates the key insight of the random algebraic constructions, that is, to provide ``very non-smooth probability distributions'' (quoted from \cite{B}).

\begin{lemma}\label{lem:dim}(\cite{B}, Lemma 5)
  For every $b$ and $s$ there exists a constant $c>0$ such the following holds: Suppose $f_1(Y),...,f_b(Y)$ are $b$ polynomials on $\F^b_q$ of degree at most $bs$, and consider the set
  $$ W = \{y \in \F^b_q : f_1(y) = \cdots = f_b(y) = 0\}.$$
  Then either $|W|< c$ or $|W|\ge q-c\sqrt q$.
\end{lemma}

We are ready to prove Theorem~\ref{thm:lower}.

\noindent {\it Proof of Theorem~\ref{thm:lower}.}
Recall the parameters defined in the beginning of this section.
We shall show that it suffices to choose the constant $c$ from Lemma \ref{lem:dim}.
Notice that $c$ only depends on $s_1,s_2,..., s_{r-1}$ and the hypergraph $H$.

As described, we choose a polynomial $f\in \sP$ uniformly at random and let $G$ be the associated $r$-uniform hypergraph $G_f$.
Let $N=q^b$ be the number of vertices in $G$, where $q$ is sufficiently large.
Following our definitions, $H$ has $v$ vertices and $e$ edges, where $\binom{e}{2}<q, \binom{v}{2}<q$ and $e< bs$.
By Lemma \ref{lem:prob}, the probability that given $v$ vertices in $G$ form a copy of $H$ is at least $1/q^e$.
Therefore, the expected number of copies of $H$ in $G$ is at least
$$\frac{1}{q^e} \cdot \binom{N}{v}=\Omega(N^{v-e/b}).$$

Let $T$ be a fixed labelled copy of $K^{(r)}_{s_1,s_2,...,s_{r-1},1}$, where we name its vertices as
$v$ and $u^i_j$'s for $1\le j\le s_i$ and $i\in [r-1]$ such that $u^i_1,\cdots, u^i_{s_i}$ are in the same part of $T$.
Fix any sequence of vertices $w^i_j$ for $1\le j\le s_i$ and $i\in [r-1]$ in $G$.
Let $W$ be the family of copies of $T$ in $G$ such that $w^i_j$ corresponds to $u^i_j$ for all $1\le j\le s_i$ and $i\in [r-1]$.
Observe that $|W|^s$ counts the number of ordered collections of $s$ copies of $T$ from $W$,
where these copies of $T$ possibly are identical.
So each such collection can be an element $L$ in
$$\mathcal{K}:=\left\{K^{(r)}_{s_1,s_2,...,s_{r-1},1},K^{(r)}_{s_1,s_2,...,s_{r-1},2},\cdots,K^{(r)}_{s_1,s_2,...,s_{r-1},s}\right\}.$$
Let $t:=s_1+\cdots+s_{r-1}$.
For given $L\in \mathcal{K}$, let $N_s(L)$ be the total number of all possible ordered collections of $s$ copies of $T\in W$ which could appear in $G$ as a copy of $L$.
So $N_s(L)=O_s\left(N^{|L|-t}\right)$.
Since the number of edges $e(L)=b\cdot (|L|-t)$ of $L$ is at most $bs$ and $q$ is sufficiently large, by Lemma \ref{lem:prob},
the probability that a potential copy $L$ appears in $G$ is $q^{-e(L)}$. Hence,
$$\E[|W|^s] =\sum_{L\in \mathcal K}N_s(L)\cdot q^{-e(L)}=\sum_{L\in \mathcal{K}}O_s\left(N^{|L|-t}\right)\cdot q^{-e(L)}
=O_s\left(\sum_{L\in \mathcal{K}}q^{b(|L|-t)}q^{-e(L)}\right)=O_s(1).$$
Note that $W$ consists of points $x\in \F^b_q$ satisfying the system of $b$ equations
$f(w^1_{j_1},w^2_{j_2},...,w^{r-1}_{j_{r-1}},x)=0$ for all $1\le j_i \le s_i$ and $i\in [r-1]$.
Because each $f(w^1_{j_1},w^2_{j_2},...,w^{r-1}_{j_{r-1}},\cdot)$ has degree at most $bs$,
by Lemma \ref{lem:dim}, either $|W|<c$ or $|W|\ge q-c\sqrt{q}\ge q/2$.
Using Markov's inequality, we then have
$$\mathbb{P}[|W|\ge c]=\mathbb{P}[|W|\ge q/2]=\mathbb{P}[|W|^s\ge (q/2)^s]\le \frac{\E[|W|^s]}{(q/2)^s}=\frac{O_s(1)}{q^s}.$$

A sequence of vertices $w^i_{j}$ for $1\le j\le s_i$ and $i\in [r-1]$ is called {\it bad},
if the corresponding set $W$ has cardinality $|W|\ge c$.
Let $B$ be the random variable counting the number of bad sequences in $G$.
Since $s=bt-b+e+1$ and $q$ is sufficiently large, it follows that
$$\E[B]\le r!N^t\cdot \frac{O_s(1)}{q^s}=O_s(q^{bt-s})=O_s(q^{b-e-1})=O_s(N^{1-\frac{e}{b}-\frac{1}{b}}).$$

We now remove a vertex from each bad sequence to form a new hypergraph $G'$.
This leaves no bad sequences in $G'$, so $G'$ is $K^{(r)}_{s_1,s_2,...,s_{r-1},p}$-free for any integer $p\ge c$.
Since each vertex is in at most $v\cdot N^{v-1}$ copies of $H$ in $G$, the total number of copies of $H$ removed is at most $B\cdot v\cdot N^{v-1}$.
Hence, the expected number of copies of $H$ in $G'$ is at least
$$\Omega(N^{v-e/b})-\E[B]\cdot v \cdot N^{v-1}=\Omega(N^{v-e/b}).$$
Therefore, for any $p\ge c$, there exists a $K^{(r)}_{s_1,s_2,...,s_{r-1},p}$-free
$r$-uniform hypergraph $G'$ with $N-o(N)$ vertices and $\Omega(N^{v-e/b})$ copies of $H$.
This completes the proof of Theorem \ref{thm:lower}.\qed

\medskip

We see that Theorem \ref{thm:lower} constructs a non-degenerate hypergraph with desired properties,
which is obtained by considering a random symmetric polynoimal $f\in \sP$.

\section{Upper bounds for degenerate hypergraphs}
In this section, we complete the proof of Theorem \ref{thm:Ksssp} by proving the following lemma.

\begin{lemma}
If $a_1<s_1\le s_2\le \cdots\le s_{r-1}$ and $a_i\le s_{i-1}$ for each $2\le i\le r$, then
$$ex(n,K^{(r)}_{a_1,a_2,...,a_r},K^{(r)}_{s_1,s_2,...,s_{r-1},s_r})\le \left(\frac{(s_r-1)^{\frac{a_1a_2\cdots a_{r}}{s_1s_2\cdots s_{r-1}}}}{\prod_{i=1}^{r} a_i!}+o(1)\right)\cdot n^{(a_1+a_2+\cdots +a_r)-\frac{a_1a_2\cdots a_r}{s_1 s_2\cdots s_{r-1}}},$$
where $o(1)$ goes to zero as $n\to \infty$.
\end{lemma}

\begin{pf}
Let $\alpha_1,...,\alpha_\ell$ be all distinct integers in $\{a_1,a_2,...,a_r\}$ and for any $j\in [\ell]$,
let $\beta_j$ be the number of all $a_i$'s which is equal to $\alpha_j$. Let $\gamma:=\beta_1!\beta_2!\cdots \beta_\ell!$.

Given an $n$-vertex $K^{(r)}_{s_1,s_2,...,s_{r-1},s_r}$-free $r$-uniform hypergraph $H$,
let $\mathcal{N}$ denote the number of ordered $r$-tuples $(A_1,...,A_r)$ of disjoint subsets of vertices such that $|A_i|=a_i$ and $A_1,...,A_r$ form a copy of $K^{(r)}_{a_1,a_2,...,a_r}$ in $H$.
For convenience, let us call such an ordered $r$-tuple $(A_1,...,A_r)$ as an {\it ordered} $K^{(r)}_{a_1,a_2,...,a_r}$.
So $\mathcal{N}$ equals the product of $\gamma$ and the number of copies of $K^{(r)}_{a_1,a_2,...,a_r}$ in $H$.

We will prove a slightly stronger statement that
$\mathcal{N}$ is at most the above right-hand side, using induction on $r$.
The base case $r=2$ follows directly from a result of Alon-Shikhelman \cite{AS} that
$$\gamma\cdot ex(n,K_{a,b},K_{s,t})\le (1+o(1))\cdot \alpha\cdot n^{a+b-ab/s},$$
where $\alpha=\frac{1}{a!(b!)^{1-a/s}}\cdot\binom{t-1}{b}^{a/s}\le \frac{(t-1)^{ab/s}}{a!b!}$.
We point out that this inequality originally states in \cite{AS} for $a\le b<s\le t$,
however the same proof actually shows that it also holds under the slightly general conditions $a\le s$ and $b<s$.
We assume that this statement holds for $(r-1)$-uniform hypergraphs.

Consider an $r$-uniform hypergraph $H$.
Let $A_1,...,A_{r-1}$ be any $r-1$ disjoint subsets of $V(H)$ of sizes $a_1,...,a_{r-1}$ respectively.
Write $\mathcal{A}= (A_1,A_2,...,A_{r-1})$ as an ordered $(r-1)$-tuple and let $n_\mathcal{A}$ be the number of vertices $w$ such that
$A_1,...,A_{r-1}, \{w\}$ induces a complete $r$-partite subhypergraph of $H$.
It then follows that
\begin{align}\label{equ:N}
\mathcal{N}\le \sum_{\mathcal{A}}\binom{n_\mathcal{A}}{a_r}\le \frac{1}{a_r!}\sum_{\mathcal{A}}n_{\mathcal{A}}^{a_r},
\end{align}
where the summations here and in what follows are over all such $\mathcal{A}$'s from $H$.
Recall the means inequality that for any $0<p\le q$, $\sum_{i=1}^m x_i^p\le m^{1-p/q}\cdot \left(\sum_{i=1}^m x_i^q \right)^{p/q}.$
Letting $m$ be the number of $(r-1)$-tuples $\mathcal{A}$ and in view of $a_r\le s_{r-1}$, we obtain from \eqref{equ:N} that
\begin{align}\label{equ:N-2}
      \mathcal{N}\le \frac{1}{a_r!}\cdot m^{1-\frac{a_r}{s_{r-1}}}\cdot \left(\sum_{\mathcal{A}}n_{\mathcal{A}}^{s_{r-1}}\right)^{a_r/s_{r-1}}\le \frac{1}{a_r!} \cdot \left(\frac{n^{a_1+\cdots+a_{r-1}}}{\prod_{i=1}^{r-1} a_i!}\right)^{1-\frac{a_r}{s_{r-1}}} \left(\sum_{\mathcal{A}}n_{\mathcal{A}}^{s_{r-1}}\right)^{\frac{a_r}{s_{r-1}}}.
\end{align}

Next we turn to prove an upper bound for $\sum_{\mathcal{A}}n_{\mathcal{A}}^{s_{r-1}}$. We claim that
\begin{align}\label{equ:nAS}
\sum_{\mathcal{A}}n_{\mathcal{A}}^{s_{r-1}}\le (c+o(1))\cdot n^{(a_1+\cdots+a_{r-1})-\frac{a_1\cdots a_{r-1}}{s_1 \cdots s_{r-2}}+s_{r-1}}, \text{ where } c=\frac{(s_r-1)^{\frac{a_1\cdots a_{r-1}}{s_1\cdots s_{r-2}}}}{\prod_{i=1}^{r-1} a_i!}.
\end{align}
For any subset $S\subseteq V(H)$ of size $s_{r-1}$, let $H_S$ be the $(r-1)$-graph on $V\backslash S$,
where $f\in E(H_S)$ if and only if $f\cup\{u\}\in E(H)$ for every $u\in S$.
Since $H$ is $K^{(r)}_{s_1,s_2,...,s_{r-1},s_r}$-free, it is clear that $H_S$ is $K^{(r-1)}_{s_1,...,s_{r-2}, s_r}$-free.
By our inductive hypothesis, the number of copies of ordered $K^{(r-1)}_{a_1,a_2,...,a_{r-1}}$ in every $H_S$ is at most
$$(c+o(1))\cdot n^{(a_1+\cdots+a_{r-1})-\frac{a_1\cdots a_{r-1}}{s_1 \cdots s_{r-2}}}.$$
Therefore we have
\begin{align}\label{equ:nAS-1}
\sum_{\mathcal{A}}\binom{n_{\mathcal{A}}}{s_{r-1}}\le (c+o(1))\cdot n^{(a_1+\cdots+a_{r-1})-\frac{a_1\cdots a_{r-1}}{s_1 \cdots s_{r-2}}} \cdot \binom{n}{s_{r-1}}.
\end{align}
Let $\mathcal{I}$ be the family consisting of all $(r-1)$-tuples $\mathcal{A}$ with $n_{\mathcal{A}}\ge \log n$, and $\mathcal{J}$ be the family of the remaining $\mathcal{A}$'s. Then by \eqref{equ:nAS-1},
\begin{equation}\label{equ:AI}
\begin{split}
\sum_{\mathcal{A}\in \mathcal{I}}n_{\mathcal{A}}^{s_{r-1}}&\le (1+o(1))\cdot s_{r-1}! \sum_{\mathcal{A}\in\mathcal{I}}\binom{n_{\mathcal{A}}}{s_{r-1}}
\le (c+o(1))\cdot s_{r-1}!\cdot n^{(a_1+\cdots+a_{r-1})-\frac{a_1\cdots a_{r-1}}{s_1 \cdots s_{r-2}}} \cdot \binom{n}{s_{r-1}}\\
&\le (c+o(1))\cdot n^{(a_1+\cdots+a_{r-1})-\frac{a_1\cdots a_{r-1}}{s_1 \cdots s_{r-2}}+s_{r-1}}.
\end{split}
\end{equation}
Since $\frac{a_1\cdots a_{r-1}}{s_1\cdots s_{r-2}}<s_{r-1},$ we have
\begin{equation}\label{equ:AJ}
\sum_{\mathcal{A}\in \mathcal{J}}n_{\mathcal{A}}^{s_{r-1}}\le n^{a_1+\cdots+a_{r-1}}\cdot (\log n)^{s_{r-1}}
=o\left(n^{(a_1+\cdots+a_{r-1})-\frac{a_1\cdots a_{r-1}}{s_1 \cdots s_{r-2}}+s_{r-1}}\right).
\end{equation}
Putting \eqref{equ:AI} and \eqref{equ:AJ} together, we derive \eqref{equ:nAS} as claimed.

Lastly, using \eqref{equ:N-2} and \eqref{equ:nAS}, it is straightforward to show that
the number $\mathcal{N}$ of copies of ordered $K^{(r)}_{a_1,a_2,...,a_r}$ in $H$ satisfies
  \begin{equation*}
    \begin{split}
     \mathcal{N}&\le \frac{1}{a_r!} \cdot \left(\frac{n^{a_1+\cdots+a_{r-1}}}{\prod_{i=1}^{r-1} a_i!}\right)^{1-\frac{a_r}{s_{r-1}}} \cdot (c+o(1))^{\frac{a_r}{s_{r-1}}}\cdot n^{(a_1+\cdots+a_{r-1})\cdot\frac{a_r}{s_{r-1}}-\frac{a_1a_2\cdots a_r}{s_1s_2\cdots s_{r-1}}+a_r}\\
    &\le \left(\frac{(s_r-1)^{\frac{a_1a_2\cdots a_{r}}{s_1s_2\cdots s_{r-1}}}}{\prod_{i=1}^{r} a_i!}+o(1)\right)\cdot n^{(a_1+a_2+\cdots +a_r)-\frac{a_1a_2\cdots a_r}{s_1 s_2\cdots s_{r-1}}}.
    \end{split}
  \end{equation*}
This completes the proof.\qed
\end{pf}

\medskip

\noindent {\bf Remark.} If $a_1=\cdots=a_r=1$, instead of requiring each $s_i\ge 2$,
we point out that one can just choose $s_1,...,s_r$ to be {\it any} positive integers and the same proof will also work.
This justifies Corollary \ref{cor:Ksssp}.

\section{Concluding remarks}
Our lower bound is obtained by the random algebraic method.
When restricting to the graph case, this provides a different construction for $ex(n,K_m,K_{s,t})$ and $ex(n,K_{a,b},K_{s,t})$,
compared to the projective norm-graph construction given in \cite{AS} (which was acquired from \cite{ARS}).
As seen in Corollaries \ref{cor:Kab-Kst} and \ref{cor:Km-Kst}, the random algebraic construction sometimes can even beyond the limit of norm-graphs.

The Ramsey numbers of complete degenerate hypergraphs also have been studied (see \cite{LM}).

It will be interesting to see if the lower bound in Theorem \ref{thm:lower} is tight when $H$ is a complete $r$-uniform hypergraph.
We are not able to get a satisfied upper bound for this case but tend to believe the answer is affirmative.

\medskip

\bigskip

\noindent {\bf Acknowledgements.} The first author is grateful to Boris Bukh for helpful discussions.
The authors would like to thank David Conlon for pointing references \cite{Er64,Mu}
and Dhruv Mubayi for informing the work of Verstra\"ete \cite{Ver} and providing several references after the first version of this draft appears on arXiv.
The authors wish to thank Verstra\"ete for illustrating his construction, which actually was done in the notes of a graduate course at UCSD in 2014.


\begin{thebibliography}{99}
\bibitem{ARS}
N. Alon, L. R\'onyai and T. Szab\'o,
Norm-graphs: variations and applications,
\textit{J. Combin. Theory Ser. B} \textbf{76} (1999), 280--290.

\bibitem{AS}
N. Alon, C. Shikhelman,
Many $T$ copies in $H$-free graphs,
\textit{J. Combin. Theory Ser. B} \textbf{121} (2016), 146--172.

\bibitem{BBK}
P. V. M. Blagojevi\'c, B. Bukh and R. Karasev,
Tur\'an numbers for $K_{s,t}$-free graphs: topological obstructions and algebraic constructions,
\textit{Israel J. Math. } \textbf{197} (2013), 199--214.


\bibitem{Boll}
B. Bollob\'as,
On complete subgraphs of different orders,
\textit{Math. Proc. Cambridge Philos. Soc.} \textbf{79} (1976), 19--24.


\bibitem{Bro}
W. Brown,
On graphs that do not contain a Thomsen graph,
\textit{Canad. Math. Bull.} \textbf{9} (1966), 281--285.


\bibitem{B}
B. Bukh,
Random algebraic construction of extremal graphs,
\textit{Bull. London Math. Soc.} \textbf{47(6)} (2015), 939--945.

\bibitem{BC}
B. Bukh and D. Conlon,
Rational exponents in extremal graph theory,
arXiv:1506.06406v1 [math.CO].

\bibitem{C}
D. Conlon,
Graphs with few paths of prescribed length between any two vertices,
\textit{Bull. London Math. Soc.}, to appear.


\bibitem{Erd}
P. Erd\H{o}s,
On the number of complete subgraphs contained in certain graphs,
\textit{Magy. Tud. Akad. Mat. Kut. Int\'ez. K\"ozl.} \textbf{7} (1962), 459--474.


\bibitem{Er64}
P. Erd\H{o}s,
On extremal problems of graphs and generalized graphs,
\textit{Israel J. Math.} \textbf{2} (1964), 183--190.


\bibitem{ERS}
P. Erd\H{o}s, A. R\'enyi and V. T. S\'os,
On a problem of graph theory,
\textit{Studia Sci. Mth. Hungar} \textbf{1} (1966), 215--235.


\bibitem{ESi}
P. Erd\H{o}s and M. Simonovits,
A limit theorem in graph theory,
\textit{Studia Sci. Mth. Hungar} \textbf{1} (1966), 51--57.

\bibitem{ESt}
P. Erd\H{o}s and A. H. Stone,
On the structure of linear graphs,
\textit{Bull. Amer. Math. Soc.} \textbf{52} (1946), 1087--1091.

\bibitem{FS}
Z. F\"uredi and M. Simonovits,
The history of degenerate (bipartite) extremal graph problems, in \textit{Erd\H{o}s Centennial},
pp. 169--264, Bolyai Soc. Math. Stud., 25, Springer, Berlin, 2013.


\bibitem{GRS}
D. Gunderson, V. R\"odl and A. Sidorenko,
Extremal problems for sets forming boolean algebras and complete partite hypergraphs,
\textit{J. Combin. Theory Ser. A} \textbf{88} (1999), 342--367.


\bibitem{GL}
E. Gy\H{o}ri and H. Li,
The maximum number of triangles in $C_{2k+1}$-free graphs,
\textit{Combin. Probab. Comput.} \textbf{21} (2012), 187--191.



\bibitem{KKM}
N. Katz, E. Krop and M. Maggioni,
Remarks on the box problem,
\textit{Math. Res. Lett.} \textbf{9} (2002), 515--519.



\bibitem{KRS}
J. Koll\'ar, L. R\'onyai and T. Szab\'o,
Norm-graphs and bipaitite Tur\'an numbers,
\textit{Combinatorica} \textbf{16} (1996), 399--406.


\bibitem{KST}
T. K\H{o}v\'ari, V. S\'os and P. Tur\'an,
On a problem of K. Zarankiewicz,
\textit{Colloquium Math.} \textbf{3} (1954), 50--57.


\bibitem{LM}
F. Lazebnik and D. Mubayi,
New lower bounds for Ramsey numbers of graphs and hypergraphs,
\textit{Adv. in Appl. Math.} \textbf{28} (2002), 544--559.


\bibitem{Mantel}
W. Mantel, Problem 28, \textit{Wiskundige Opgaven} \textbf{10} (1907), 60--61.


\bibitem{Mu}
D. Mubayi,
Some exact results and new asymptotics for hypergraph Tur\'an numbers,
\textit{Combin. Probab. Comput.} \textbf{11} (2002), 299--309.


\bibitem{Sim}
M. Simonovits,
A method for solving extremal problems in graph theory, stability problems,
\textit{Theory of Graphs (Proc. Colloq., Tihany, 1966)}, 279--319.


\bibitem{Tur}
P. Tur\'an,
On an extremal problem in graph theory (in Hungarian),
\textit{Math. Fiz. Lapok} \textbf{48} (1941), 436--452.

\bibitem{Ver}
J. Verstra\"ete,
Private communication.


\end{thebibliography}
\end{document}